\documentclass[11pt]{article}

\usepackage{epsfig}
\usepackage{amssymb}
\usepackage{amsthm}

\newtheorem{remark}{Remark}

\usepackage{mathtools}
\usepackage{relsize}
\usepackage{xcolor}

\def\Frac#1#2{\frac{\displaystyle{#1}}{\displaystyle{#2}}}
\def\protectbold#1{\protect{\boldmath{$#1$}}}

\def\eqref#1{(\ref{#1})}
\def\tfrac#1#2{{{\lower.6ex
\hbox{$\scriptstyle#1$}}\over
{\raise.7ex
\hbox{$\scriptstyle#2$}}}}
\def\sign{{\rm sign}}
\def\erfc{{\rm erfc}}

\def\phase{{\rm ph}}
\def\bigO{{\cal O}}
\def\RR{{\mathbb R}}
\def\dsp#1{\displaystyle#1}

\def\FG#1#2#3#4{
{}_2F_1\left(
\begin{array}{c}
\begin{array}{c}\hskip-10pt#1,#2\end{array}\\
\begin{array}{c}\hskip-10pt #3\end{array}
\end{array}
\hskip-8pt;\,#4
\right)}

\numberwithin{equation}{section}

\begin{document}

 \title{On the computation and inversion of the cumulative noncentral beta distribution function}

\author{Amparo Gil\\
Departamento de Matem\'atica Aplicada y CC. de la Computaci\'on.\\
ETSI Caminos. Universidad de Cantabria. 39005-Santander, Spain.\\
\and
Javier Segura\\
        Departamento de Matem\'aticas, Estad\'{\i}stica y 
        Computaci\'on,\\
        Univ. de Cantabria, 39005 Santander, Spain.\\   
\and
Nico M. Temme\\
IAA, 1825 BD 25, Alkmaar, The Netherlands\footnote{Former address: Centrum Wiskunde \& Informatica (CWI), 
        Science Park 123, 1098 XG Amsterdam,  The Netherlands}.\\ \\
{ \small
  e-mail: {\tt
    amparo.gil@unican.es,
    javier.segura@unican.es, 
    nico.temme@cwi.nl}}
    }

\date{\ }

\maketitle
\begin{abstract}
The computation and inversion of the noncentral beta distribution $B_{p,q}(x,y)$ 
(or the noncentral $F$-distribution, a particular case of  $B_{p,q}(x,y)$)
play an important role
in different applications. In this paper we study the stability of recursions satisfied by
 $B_{p,q}(x,y)$  and its complementary function and describe asymptotic expansions useful for
computing the function when the parameters are large. We also
consider the inversion problem of finding $x$ or $y$ when a value of $B_{p,q}(x,y)$ is given.
We provide approximations to $x$ and $y$ which can be used as starting values
of methods for solving nonlinear equations (such as Newton) if higher accuracy is needed. 
 
\end{abstract}

\vskip 0.5cm \noindent
{\small
\par\noindent
Keywords \& Phrases:
Noncentral beta distribution; asymptotic expansions; numerical computation; inversion.
}

\section{Introduction}\label{sec:intro}

The cumulative noncentral  beta distribution can be defined as\footnote{Different notations are used in the literature; for this we refer to \S\ref{sec:notations}. }
\begin{equation}
\label{defi}
B_{p,q}(x,y)=e^{-x/2} \displaystyle\sum_{j=0}^{\infty} \Frac{1}{j!}\left(\Frac{x}{2}\right)^j I_y (p+j,q),
\end{equation}
where $I_y (p,q)$ is the central beta distribution
\begin{equation}
\label{cent}
I_y (p,q)=\Frac{1}{B(p,q)}\displaystyle\int_0^y t^{p-1} (1-t)^{q-1} dt,\quad 
B(p,q)=\Frac{\Gamma (p)\Gamma (q)}{\Gamma (p+q)}.
\end{equation}
 These functions satisfy the recurrence relation
\begin{equation}
\label{recIy}
(p+j)f_{j+1}=\left(p+j+(p+j+q-1)y\right)f_j-(p+j+q-1)yf_{j-1},
\end{equation}
where $f_j=I_y (p+j,q)$. The recursion is not stable in the forward direction, but we can use a backward recursion scheme (see \cite[\S4.6]{Gil:2007:NSF}).

The complementary function is defined by 
\begin{equation}\label{Bbar}
\bar{B}_{p,q}(x,y)=1-B_{p,q}(x,y),
\end{equation}
and using $I_{y}(p,q)=1-I_{1-y}(q,p)$ we
can write the series
\begin{equation}\label{Bbarexpand}
\bar{B}_{p,q}(x,y)=e^{-x/2}\displaystyle\sum_{j=0}^{\infty} \Frac{1}{j!}\left(\Frac{x}{2}\right)^j I_{1-y}(q,p+j).
\end{equation}
This is related with the so-called beta distribution of type II, $J_x (p,q,\lambda)$ \cite{Chat:1995:ANO}, by 
$J_y (p,q,x)=\bar{B}_{q,p}(1-y,x)$.

Using (\ref{cent}) in (\ref{defi}) we get
\begin{equation}
B_{p,q}(x,y)=\Frac{e^{-x/2}}{B(p,q)}\displaystyle\int_0^y t^{p-1} (1-t)^{q-1} M\left(p+q,p,\tfrac12xt\right)\,dt,
\end{equation}
where $M(a,b,z)$ is the confluent hypergeometric function or Kummer function; see \cite[Ch.~13]{Daalhuis:2010:CHF}. Similarly for $\bar{B}_{p,q}(x,y)$ but with the integration going from $y$ to $1$. 

Therefore the probability density function can be written
\begin{equation}
\label{pary}
\Frac{\partial B_{p,q}(x,y)}{\partial y}=\Frac{e^{-x/2}}{B(p,q)} y^{p-1} (1-y)^{q-1} M\left(p+q,p,\tfrac12xy\right)>0.
\end{equation}

The derivative with respect to $x$ can be computed using (\ref{defi}). We get
\begin{equation}
\Frac{\partial B_{p,q}(x,y)}{\partial x}=\tfrac{1}{2}\left(B_{p+1,q}(x,y)-B_{p,q}(x,y)\right),
\end{equation}
which, using the recurrence relation (\ref{rede}) that we later obtain, can be written
\begin{equation}
\label{parx}
\Frac{\partial B_{p,q}(x,y)}{\partial x}=-\Frac{e^{-x/2}}{2 p B(p,q)} y^{p} (1-y)^{q} M\left(p+q,p+1,\tfrac12xy\right)<0.
\end{equation}
By using the recurrence relations for the Kummer functions \cite[\S13.3(i)]{Daalhuis:2010:CHF} it is possible to check, that, as
expected, taking the derivative of (\ref{pary}) with respect to $x$ gives the same value as the derivative of (\ref{parx})
with respect to $y$.

Particular values are (some are evident and the rest become evident in a later analysis)
\begin{equation}
\begin{array}{ll}
B_{p,q}(x,0)=0, & B_{p,q}(x,1)=1,\\[6pt]
B_{p,q}(0,y)=I_y (p,q), & B_{p,q}(+\infty,y)=0,\\[6pt]
B_{+\infty,q}(x,y)=0, & B_{p,+\infty}(x,y)=1,  
\end{array}
\end{equation}
and, as corresponds to a cumulative distribution function, we have $0\le B_{p,q}(x,y)\le 1$ for $p,q>0$ and $x,y\ge 0$.

In this paper, we give useful expressions 
for computing $B_{p,q}(x,y)$ and its complementary function: inhomogeneous and three and four-term homogeneous
recurrence relations; asymptotic expansions.
For the three and four-term recurrence relations, their stability is also analyzed. Numerical examples illustrate the accuracy of  
the expressions.
We also consider the inversion problem of finding $x$ or $y$ when a value of $B_{p,q}(x,y)$ is given.
We provide approximations to $x$ and $y$ which can be used as starting values
of methods for solving nonlinear equations (such as Newton) if higher accuracy is needed.

\section{Other notations in the literature}\label{sec:notations}
We have chosen for the definition and notation in \eqref{defi}, knowing that several other definitions and notations exist. We like to separate the parameters $p,q$ from $x,y$. The first set can be used for recursions, the second set for differentiation. 

First of all, the noncentral beta distribution and the noncentral-$F$ distribution are the same, up to notation of the parameters. The cumulative distribution function $F(w,\nu_1,\nu_2,\lambda)$ of the noncentral $F$-distribution, with $\nu_1,\nu_2$ degrees of freedom and noncentrality parameter $\lambda$ and the noncentral beta distribution are related as follows:
\begin{equation}\label{eq:rel01}
F(w,\nu_1,\nu_2,\lambda)=B_{p,q}(\lambda,x),
\end{equation}
where
\begin{equation}\label{eq:rel02}
\nu_1=2p,\quad \nu_2=2q,\quad x=\frac{\nu_1 w}{\nu_1w + \nu_2}.
\end{equation}

In \cite{Robertson:1976:CNF} and \cite{Helstrom:1985:ENF} a function $q(\lambda,\omega;2a,2b)$ is defined by
\begin{equation}\label{eq:rel03}
1-q(\lambda,\omega;2a,2b)=e^{-\lambda}x^{a}(1-x)^b\sum_{n=0}^\infty \frac{\Gamma(a+b+n)}{\Gamma(b)\Gamma(a+1+n)}x^n\sum_{j=0}^n\frac{\lambda^j}{j!},
\end{equation}
where $ x=\omega/(\omega+1)$. In our notation we have
\begin{equation}\label{eq:rel04}
1-q(\lambda,\omega;2a,2b)=B_{a,b}(2\lambda,x).
\end{equation}
For a proof we refer to the Appendix. It follows that (see \eqref{Bbarexpand})
\begin{equation}\label{eq:rel05}
q(\lambda,\omega;2a,2b)=\bar{B}_{a,b}(2\lambda,x).
\end{equation}

In \cite{Chat:1995:ANO},
\cite{Baharev:2008:CNF} and \cite{Baharev:2017:CNF} the definition of the noncentral beta function is\begin{equation}\label{eq:rel06}
I_x(a,b;\lambda)=e^{-\frac12\lambda}\sum_{j=0}^\infty\frac{1}{j!}\left(\frac{x}{2}\right)^jI_x(a+j,b),
\end{equation}
and the relation with our definition is
\begin{equation}\label{eq:rel07}
I_x(a,b;\lambda)=B_{a,b}(\lambda,x).
\end{equation}

\section{Recurrence relations}\label{sec:recrel}

We obtain two-term inhomogeneous recurrence relations for the beta distributions, as well as three and
 four-term homogeneous recurrence relations, and we discuss the stability of these recursive processes.

\subsection{First order recurrence relations}\label{sec:firstR}

Using the recurrence relations for the central distribution given in \cite[8.17(iv)]{Paris:2010:IGA} together with the
definition (\ref{defi}) it is easy to get analogous relations for the noncentral case.

Using \cite[8.17.20]{Paris:2010:IGA} we obtain
\begin{equation}
\label{rede}
B_{p,q}(x,y)=B_{p+1,q}(x,y)+\Frac{e^{-x/2}}{pB(p,q)}y^p (1-y)^q M \left(p+q,p+1,\tfrac12xy\right),
\end{equation}
from \cite[8.17.21]{Paris:2010:IGA}
\begin{equation}
\label{redb}
B_{p,q}(x,y)=B_{p,q+1}(x,y)-\Frac{e^{-x/2}}{qB(p,q)}y^p (1-y)^q M \left(p+q,p,\tfrac12xy\right),
\end{equation}
and from \cite[8.17.19]{Paris:2010:IGA}
\begin{equation}
B_{p,q}(x,y)=B_{p-1,q+1}(x,y)-\Frac{e^{-x/2}}{qB(p,q)}y^{p-1} (1-y)^q M \left(p+q,p,\tfrac12xy\right),
\end{equation}

For the complementary function, $\bar{B}_{p,q}(x,y)=1-B_{p,q}(x,y)$ the same recursions hold but with 
the sign of the inhomogeneous term reversed, for instance, the first recursion reads
\begin{equation}
\label{rede2}
\bar{B}_{p,q}(x,y)=\bar{B}_{p+1,q}(x,y)-\Frac{e^{-x/2}}{pB(p,q)}y^p (1-y)^q M \left(p+q,p+1,\tfrac12xy\right).
\end{equation}

In view of the recursions, it appears that $B_{p,q}(x,y)$ ($\bar{B}_{p,q}(x,y)$) can be computed in the direction of decreasing (increasing) 
$p$ and in the direction of increasing (decreasing) $q$. In turn, this means that $B_{p,q}(x,y)$ is minimal with respect to $p$ and $\bar{B}_{p,q}(x,y)$
is minimal with respect to $q$.

\subsection{Homogeneous three-term recurrence relations}\label{sec:hom3}

Next we provide equivalent three-term recurrence relations which display similar stability as the two-term
relations but which have some advantages, as we next discuss. 

For instance, using the recurrence relation (\ref{rede}) we can write 
$(B_{p+1,q}-B_{p,q})/(B_{p,q}-B_{p-1,q})$ in terms of a ratio of Kummer functions, and this leads to

\begin{equation}
\label{TTRR1}
\begin{array}{l}
y_{p+1}-(1+c_{p,q})y_{p}+c_{p,q}y_{p-1}=0,\\
\\
c_{p,q}=\Frac{p+q-1}{p}y\Frac{M(p+q,p+1,xy/2)}{M(p+q-1,p,xy/2)},
\end{array}
\end{equation}
where $y_{p+j}=B_{p+j,q}$. The three-term recurrence is also satisfied by $y_{p+j}=\bar{B}(p+j,q)$.

The coefficient $c_{p,q}$ can be computed in terms of a continued fraction because the Kummer function is minimal for the recurrence relation in the
$(+,+)$ direction (see \cite{Segura:2008:NSS}); in addition, using the behavior predicted by the 
Perron-Kreuser theorem (see \cite{Kreuser:1914:UDV} or \cite{Cash:1980:ANO} for a more 
recent account of this useful theorem) 
for this recurrence relation, we
obtain that $c_{p,q}(x,y)\rightarrow y$ as $p\rightarrow +\infty$. The application of the Perron-Kreuser theorem to (\ref{TTRR1}), 
gives us information on the behavior as $p\rightarrow +\infty$ and we conclude that as $p\rightarrow +\infty$,
\begin{equation}
\label{compoasp}
\Frac{B_{p+1,q}}{B_{p,q}}\sim c_{p,q}\sim y, \quad \Frac{\bar{B}_{p+1,q}}{\bar{B}_{p,q}}\sim 1
\end{equation}
and, as announced,  $B_{p,q}$ is minimal while $\bar{B}_{p,q}$ is dominant. Because $B_{p,q}$ is minimal, Pincherle's theorem
gives the following continued fraction representation
$$
\Frac{B_{p+1,q}}{B_{p,q}}=\Frac{c_{p,q}}{1+c_{p,q}-}\Frac{c_{p+1,q}}{1+c_{p+1,q}-}\ldots,
$$
where, as commented before, the coefficients $c_{p,q}$ can also be computed using a continued fraction representation.

Proceeding similarly with the recurrence relation (\ref{redb}) we arrive to the three-term recursion

\begin{equation}
\label{TTRR2}
\begin{array}{l}
 y_{q+1}-(1+\bar{c}_{p,q}) y_{q}+\bar{c}_{p,q} y_{q-1}=0,\\
\\
\bar{c}_{p,q}=\Frac{p+q-1}{q}(1-y)\Frac{M(p+q,p,xy/2)}{M(p+q-1,p,xy/2)}
\end{array}
\end{equation}
where $y_{q+j}=B_{p,q+j}$ or $y_{q+j}=\bar{B}_{p,q+j}$.

We conclude that as $q\rightarrow +\infty$,
\begin{equation}
\label{compoasq}
\Frac{\bar{B}_{p,q+1}}{\bar{B}_{p,q}}\sim \bar{c}_{p,q}\sim 1-y, \quad\Frac{B_{p,q+1}}{B_{p,q}}\sim 1
\end{equation}
and, as announced,  $\bar{B}_{p,q}$ is minimal while $B_{p,q}$ is dominant. As before, there is a continued fraction representation
for $\bar{B}_{p,q+1}/\bar{B}_{p,q}$.

These recurrence relations can be useful for computing both the beta and the complementary beta distributions,
provided the recurrences are applied in their numerically stable directions: decreasing $p$ and increasing 
$q$ for $B_{p,q}$ and the opposite directions for $\bar{B}_{p,q}$. There is, however, the added difficulty
that the coefficients $c_{p,q}$ or $\bar{c}_{p,q}$ are not simple and should be evaluated by some means
 (for instance, with the associated continued fraction). These coefficients are, on the other hand, 
easier to compute and less prone to possible overflows/underflows than the inhomogeneous terms in 
(\ref{rede}) and (\ref{redb}).

\subsection{Homogeneous four-term recurrence relations}\label{sec:hom4}

Next, we provide four-term recurrence relations for $B_{p,q}$ and $\bar{B}_{p,q}$. These recurrences have 
some advantages, but also some disadvantages. The main advantage is that the coefficients are simple. 
The disadvantages are, on one hand that three function values are needed to start the recurrence, and on the
other hand that their applicability is more restricted for reasons of stability; as we next discuss, apart
from $B_{p,q}$ and $\bar{B}_{p,q}$, these recurrences introduce an independent third independent solution
which is spurious and that in some cases compromises the stability. 

For obtaining these recurrences we note that from \eqref{rede} (similarly \eqref{redb}), we can write 
the Kummer function on the right as the difference of $B_{p+1,q}$ and $B_{p,q}$. Similarly, we can do the
same for two additional consecutive Kummer functions by shifting $p$. Then, relating these three consecutive
Kummer functions by its corresponding three-term recurrence relation \cite[\S13.3(i)]{Daalhuis:2010:CHF},
we can find a recurrence for four consecutive values $B_{p+j,q}$, $j=0,1,2,3$. We get:

\begin{equation}\label{rerecurra}
c_3y_{p+3}+c_2y_{p+2}+c_1y_{p+1}+c_0y_{p}=0,
\end{equation}
where
\begin{equation}\label{c0c1c2c3} 
\begin{array}{@{}r@{\;}c@{\;}l@{}}
c_0&=& (p+q)y^2,\\[8pt]
c_1&=&-y\left(p+1-\frac12xy+y(p+q)\right),\\[8pt]
c_2&=& y\left(p+1-\frac12xy\right)-\frac12xy,\\[8pt]
c_3&=&\frac12xy,
\end{array}\end{equation}
and $y_{p+j}=B_{p+j,q}(x,y)$ (or $y_{p+j}=\bar{B}_{p+j,q}(x,y)$).

Considering now the Perron-Kreuser theorem, we observe that the 
dominant solution is such that
\begin{equation}
\label{domi4}
\displaystyle\lim_{p\rightarrow \infty} \Frac{1}{p}\Frac{y_{p+1}}{y_p}=-\Frac{2}{x}.
\end{equation}
This solution is spurious because it is not any of the solutions we are interested in. Therefore, this 
recurrence should never be used in the direction of increasing $p$, neither for computing 
$B_{p,q}(x,y)$ or $\bar{B}_{p,q}(x,y)$. Then, we have a second subspace of solutions, of 
dimension two, with all the solutions such that $\limsup |y_n|^{1/n}$ is finite. The characteristic
equation for this subspace of solutions is $y\lambda^2 -y(1+y)\lambda+y^2=0$, with solutions $\lambda =1$
and $\lambda=y$. This corresponds to $y_{p}=B_{p,q}(x,y)$ for $\lambda=y$ and to 
$y_p=\bar{B}_{p,q}(x,y)$ (or $y_p=1$) for $\lambda =1$ (see \eqref{compoasp}). Because $0<y<1$, 
$B_{p,q}(x,y)$ is minimal as $p$ increases, while $\bar{B}_{p,q}(x,y)$ is dominant over $B_{p,q}(x,y)$,
 but it is not a dominant solution of the recurrence, which increases much faster according to 
(\ref{domi4}). We therefore conclude that this recurrence can be only used for computing 
$B_{p,q}(x,y)$ and only in the backward direction (decreasing $p$); the recurrence 
should not be used for $\bar{B}_{p,q}(x,y)$ neither in the forward nor in the backward direction.

As for the recurrence over $q$, we have
\begin{equation}\label{rerecurrb}
d_3y_{q+3}+d_2y_{q+2}+d_1y_{q+1}+d_0y_{q}=0,
\end{equation}
where
\begin{equation}\label{d0d1d2d3} 
\begin{array}{@{}r@{\;}c@{\;}l@{}}
d_0&=& -(p+q)(1-y)^2,\\[8pt]
d_1&=&(1-y)\left((p+q)(1-y)+p+2q+2+\frac12xy\right),\\[8pt]
d_2&=& -(1-y)\left(p+2q+2+\frac12xy\right)-q-2,\\[8pt]
d_3&=&q+2,
\end{array}\end{equation}
where we can have $y_{q+j}=B_{p,q+j}(x,y)$ or $y_{q+j}=\bar{B}_{p,q+j}(x,y)$.

The coefficients are all linear in $q$ and this implies that the characteristic equation is of third
degree, namely: $(\lambda-1)(\lambda-(1-y))^2 =0$. We have the simple root $\lambda=1$ which corresponds
to $B_{p,q}$ because of (\ref{compoasq}); therefore this is a dominant solution and it can be computed 
by forward recursion (increasing $q$). With respect to the second (double) characteristic root $\lambda=1-y$, 
the Perron-Kreuser is not so conclusive, and the information we have is that there is a subspace
of solutions of dimension two associated to this root, and that all the solutions in this subset satisfy
$$
\displaystyle\limsup_{n\rightarrow \infty}|y_q|^{1/q}=|1-y|,
$$
but we cannot conclude from here that $\bar{B}_{p,q}$ is minimal (and, in fact, it does not appear to be).
Therefore, we conclude that the recurrence can only be safely used for $B_{p,q}$ and only in the forward
direction (increasing $q$).

\subsection{Computational examples}

Examples of the numerical performance of the homogeneous three-term recurrence relations  \eqref{TTRR1} and 
 \eqref{TTRR2}, are shown Table~\ref{table1}.  The initial function values needed for starting the recurrences
and the final function values used for testing, have been computed with Maple with a large number of digits.
The recurrences have been implemented with $Digits:=16$ (double precision). 
As can be seen, the recursions are stable when applied in the correct
direction. We have checked that similar results are obtained when using the first order recurrences discussed in Section
\ref{sec:firstR}. 

\begin{table}
$$
\begin{array}{llll}
\hline
 \mbox{TTRR} &  \mbox{Function/Direction}    & \mbox{IPVs/Accuracy} & \mbox{FPVs/Accuracy} \\ 
  \hline
\eqref{TTRR1} & B_{p,q}(50,0.4)  &  p=300, q=200 & p=50, q=200   \\
                     &    \mbox{Backward}                   &  \sim 5\, 10^{-16} & 1.9\,10^{-14} \\
\hline
\eqref{TTRR1} & \bar{B}_{p,q}(50,0.4)  &  p=30, q=200 & p=280, q=200   \\
                     &    \mbox{Forward}                   &  \sim 5\, 10^{-16} & 7\,10^{-15} \\
\hline
\eqref{TTRR2} & B_{p,q}(50,0.4)  &  p=30, q=20 & p=30, q=270   \\
                     &    \mbox{Forward}                   &  \sim 5\, 10^{-16} & 4.6\,10^{-14} \\
\hline
\eqref{TTRR2} & \bar{B}_{p,q}(50,0.4)  &  p=30, q=300 & p=30, q=50   \\
                     &    \mbox{Backward}                   &  \sim 5\, 10^{-16} & 1.5\,10^{-14} \\
\end{array}
$$
{\footnotesize {\bf Table~1}. Examples of computations
using the homogeneous three-term recurrence relations \eqref{TTRR1} and \eqref{TTRR2}.
 In the second column,  the function ($B_{p,q}(x,y)$ or $\bar{B}_{p,q}(x,y)$)
and the direction of application of the TTRR (forward or backward) are specified.
The third column shows the initial values (IPVs) of the parameters $p$ and $q$ and the accuracy
of the function values.
The fourth column shows the final parameter values (FPVs) after applying
the TTRR and the accuracy of the function values.
 }
\label{table1}
\end{table}

As an example of the performance of the four-term recurrence relation
for computing $B_{p,q}(x,y)$, we consider \eqref{rerecurra} for $p=1000$ and $q=1200,\,x=10,\,y=0.2$.
As we will later show in the paper, using very few terms in the asymptotic expansion (\ref{eq:ras10}) we can compute $B_{p,q}(x,y)$ with
an accuracy close to single precision $\sim 10^{-8}$ for these large values of the parameters $p$ and $q$. 
For using the recurrence relation  (\ref{rerecurra})
three starting values are needed:  $B_{p,1200}(10,0.2)$ for $p=1003,\,1002,\,1001$. We evaluate these functions with
an accuracy of $3.4\,10^{-8}$ for the three values of $p$ using just $3$ terms in the asymptotic expansion (\ref{eq:ras10}). Then, we
use (\ref{eq:ras10}) backwards down to $p=200$. When checking the accuracy of the last function value computed 
($B_{200,1200}(10,0.2)$) we find that this is exactly $3.4\,10^{-8}$, i.e. there isn't any loss of accuracy when using  (\ref{rerecurra}).
This shows the numerical stability of the recursion. Furthermore, because the solution is minimal, the ratio
of consecutive functions tends to converge to the exact solution; in our example, which we compute with 
$16$ digits, all digits of $B_{200,1200}(10,0.2)/B_{199,1200}(10,0.2)$ are correct.

As mentioned before, for computing the function $\bar{B}_{p,q}(x,y)$ the use of  the four-term recurrence relations should always 
be avoided.
As an example, let us consider \eqref{rerecurra} for $p=300$, $q=200$, $x=10$, $y=0.2$.  The three starting values needed to apply the
recursion forward ($\bar{B}_{300,200}(10,0.2)$,
  $\bar{B}_{301,200}(10,0.2)$,  $\bar{B}_{302,200}(10,0.2)$) have been computed with Maple using a large number of digits.
When implemented in double precision, we
obtain a loss of accuracy of four significant digits in the function value when the recurrence relation \eqref{rerecurra} is used to compute  $\bar{B}_{307,200}(10,0.2)$ (just five steps of the recursive process has been applied.) 

\section{Expansions in terms of Kummer functions}\label{sec:expkum}

The first-order recursions, when iterated, lead to series in terms of Kummer functions that can be useful
for computing the beta distributions, particularly for small $y$.

Observe that the recursion in (\ref{rede}) has the form  $B_{p+1,q}(x,y)=B_{p,q}(x,y)+h_p$, and that, hence,
$B_{p,q}(x,y)=B_{p+N+1,q}(x,y)+h_p+h_{p+1}+\ldots+h_{p+N+1}$, $N=0,1,2,\ldots$. This gives
\begin{equation}\label{Nseries} 
\begin{array}{@{}r@{\;}c@{\;}l@{}}
&&B_{p,q}(x,y)=B_{p+N+1,q}(x,y)+\\[8pt]
&&\quad\quad\Frac{e^{-x/2}y^p (1-y)^q}{pB(p,q)}\displaystyle\sum_{j=0}^N y^j \Frac{(p+q)_j}{(p+1)_j} M\left(p+q+j,p+1+j,\tfrac12xy\right).
\end{array}
\end{equation}
Taking into account (\ref{compoasp}) it
is clear that $B_{p+N+1,q}\rightarrow 0$ as $N\rightarrow +\infty$. Therefore
\begin{equation}\label{Infseries1} 
B_{p,q}(x,y)=\Frac{e^{-\frac12x}y^p (1-y)^q}{pB(p,q)}\displaystyle\sum_{j=0}^{\infty} y^j  \Frac{(p+q)_j}{(p+1)_j} M\left(p+q+j,p+1+j,\tfrac12xy\right).
\end{equation}
Using $M(a,b,z)=e^zM(b-a,b,-z)$, we obtain
\begin{equation}\label{Infseries2} 
B_{p,q}(x,y)=\Frac{e^{\frac12x(y-1)}y^p (1-y)^q}{pB(p,q)}\displaystyle\sum_{j=0}^{\infty} y^j  \Frac{(p+q)_j}{(p+1)_j} M\left(1-q,p+1+j,-\tfrac12xy\right).
\end{equation}
Observing that $M(a,b+\lambda,z)=\bigO(1)$ as $\lambda\to+\infty$, we conclude again that the series is convergent if $0\le y<1$.

By using  (see \cite[\S13.3(ii)]{Daalhuis:2010:CHF})
\begin{equation}\label{Infseries3}
\Frac{d^j M(a+b,a+1,z)}{dz^j}=\Frac{(a+b)_j}{(a+1)_j} M\left(a+b+j,a+1+j,z\right),
\end{equation}
it follows that we can write
\begin{equation}\label{Kumdiff} 
B_{p,q}(x,y)=\Frac{e^{-x/2}y^p (1-y)^q}{pB(p,q)}\displaystyle\sum_{j=0}^\infty y^j \Frac{d^j M(p+q,p+1,z)}{dz^j},
\quad z=\tfrac12xy.
\end{equation}

For the complementary function we obtain the series expansion
\begin{equation}\label{Infseries4}
\bar{B}_{p,q}(x,y)=\Frac{e^{-x/2}y^p (1-y)^q}{qB(p,q)} \displaystyle\sum_{j=0}^{\infty}
(1-y)^j\Frac{(p+q)_j}{(q)_j} M\left(p+q+j,p,\tfrac12xy\right).
\end{equation}

\section{Contour integral representation}\label{sec:contour}
For deriving asymptotic expansions it is convenient to have integral representations in terms of elementary functions. We use (see \cite[\S13.4(ii)]{Daalhuis:2010:CHF}) the contour integral
\begin{equation}\label{intKum} 
M(a,b,z)=\frac{\Gamma(b)\Gamma(1+a-b)}{\Gamma(a)}\frac{1}{2\pi i}\int_{0}^{(1+)}e^{zt}t^{a-1}{(t-1)^{b-a-1}}\,dt,
\end{equation}
where $a>0$,  $b-a\not= 1,2,3,\ldots$, and the contour starts at the origin, encircles the point $t=1$ anti-clockwise, and returns to the origin. The multivalued factors
$t^{a-1}$  and $(t-1)^{b-a-1}$ assume their principal values.

When we use the expansion in \eqref{Kumdiff} and the representation in \eqref{intKum} we obtain, after interchanging the order of summation and integration, the contour integral
\begin{equation}\label{intB1} 
B_{p,q}(x,y)=\frac{e^{-x/2}y^p (1-y)^q}{2\pi i}
\int_{0}^{(1+)}e^{\frac12xyt}t^{p+q-1}{(t-1)^{-q}}\,\frac{dt}{1-yt},
\end{equation}
where we have assumed that $\vert yt\vert <1$. It follows that the point $t=1/y$ should be outside the contour, that is, the contour cuts the positive real axis between $1$ and $1/y$.

When we extend the contour in order to include the point $t=1/y$ inside the contour, calculate the residue (which equals 1), we obtain
\begin{equation}\label{intB2} 
B_{p,q}(x,y)=1+\frac{e^{-x/2}y^p (1-y)^q}{2\pi i}
\int_{0}^{(1+,1/y+)}e^{\frac12xyt}t^{p+q-1}{(t-1)^{-q}}\,\frac{dt}{1-yt},
\end{equation}
where the contour cuts the real axis on the right of the point $t=1/y$. It follows that we have for the complementary function the representation
\begin{equation}\label{intB3} 
\bar{B}_{p,q}(x,y)=\frac{e^{-x/2}y^p (1-y)^q}{2\pi i}
\int_{0}^{(1+,1/y+)}e^{\frac12xyt}t^{p+q-1}{(t-1)^{-q}}\,\frac{dt}{yt-1}.
\end{equation}

\section{Asymptotic expansions}\label{sec:asexp}
We introduce the notations
\begin{equation}\label{eq:as01} 
 p=r\cos^2\theta, \quad q=r \sin^2\theta, \quad r=p+q,\quad 0\le\theta\le\tfrac12\pi, \quad z=\tfrac12xy=r\xi. 
\end{equation}
We derive an asymptotic expansion for large values of $z$ with $p$ and $q$ bounded, and one with large values of $r$ with $\xi$ bounded.

\subsection{Asymptotic expansion for large values of \protectbold{z}}\label{sect:zas}
For large values of $z$ we use the integral representation in \eqref{intB1} and apply Watson's lemma for contour integrals, see \cite[\S2.2]{Temme:2015:AMI}. In the lemma an integral is considered of the form
\begin{equation}\label{eq:zas01} 
G_\lambda(z)=\frac{1}{2\pi i}\int_{-\infty}^{(0+)} s^{\lambda-1} q(s) e^{zs}\,ds,
\end{equation}
where  $\lambda$ is a real or complex constant and the path runs around the branch cut of the multivalued function $s^{\lambda-1}$  from $-\infty$ to $0$. Along the lower side of the negative real axis we have $\phase\,s=-\pi$, along the upper side $\phase\,s=\pi$. It is assumed that  $q(s)$ is analytic at the origin and inside the contour. The expansion $\dsp{q(s)= \sum_{n=0}^\infty a_n s^n} $ gives 
 the asymptotic expansion
\begin{equation}\label{eq:zas02}
G(z)\sim  \sum_{n=0}^\infty \frac{1}{\Gamma(1-\lambda-n)}\,\frac{a_n}{z^{n+\lambda}},\quad z\to\infty.
\end{equation}
In this result all fractional powers have their principal values. 

The integral in \eqref{intB1} is not of the form in \eqref{eq:zas01}, we have a finite contour starting at $t=0$, but it is clear that the dominant exponential contributions $e^{zt}$ are coming from $t>1$. We expand
\begin{equation}\label{eq:zas03}
t^{p+q-1}=\sum_{n=0}^\infty a_n (t-1)^n,\quad a_n=(-1)^n\frac{(1-p-q)_n}{n!},
\end{equation}
and 
\begin{equation}\label{eq:zas04}
\frac{1}{1-yt}=\frac{1}{1-y}\sum_{n=0}^\infty b_n (t-1)^n,\quad b_n=\left(\frac{y}{1-y}\right)^n.
\end{equation}
This gives
\begin{equation}\label{eq:zas05}
\frac{t^{p+q-1}}{1-yt}=\frac{1}{1-y}\sum_{n=0}^\infty c_n (t-1)^n,\quad c_n=\sum_{m=0}^n a_mb_{n-m}.
\end{equation}

When we use this expansion for the integral in \eqref{intB1}, we obtain
\begin{equation}\label{eq:zas06} 
B_{p,q}(x,y)\sim \frac{e^{-(1-y)x/2}y^p (1-y)^{q-1} z^{q-1}}{\Gamma(q)}\sum_{n=0}^\infty(-1)^n(1-q)_n\frac{c_n}{z^n},\quad z\to\infty,
\end{equation}
where $z$ is defined in \eqref{eq:as01}. 
 This expansion can be used for $p$ and $q$ of order ${\cal{O}}(1)$ and for …
$0 < y \le 1-\delta$, where $\delta$ is a fixed positive small number.

\begin{remark}\label{rem:01}
{\rm
The expansion in \eqref{eq:zas06} becomes finite and exact when $q$ is a positive integer. In that case it corresponds with the finite expansion given in \cite[Eq.~(11)]{Baharev:2008:CNF}. We have verified and confirmed this relation, although the representation of the two results is quite different, not only because of a different notation.
}
\end{remark}

\subsection{Asymptotic expansion for large values of \protectbold{r}}\label{sect:ras}
We write the integral in \eqref{intB1} in the form
\begin{equation}\label{eq:ras01} 
B_{p,q}(x,y)=\frac{e^{-x/2}y^p (1-y)^q}{2\pi i}
\int_{0}^{(1+)}e^{r\phi(t)}\,\frac{dt}{t(1-yt)},
\end{equation}
where
\begin{equation}\label{eq:ras02} 
\phi(t)=\ln t-\sin^2\theta\ln(t-1)+\xi t, \quad \phi^\prime(t)=\frac{\xi t^2+\left(\cos^2\theta-\xi\right)t-1}{t(t-1)}.
\end{equation}
There are two saddle points, one is negative, the other one is given by
\begin{equation}\label{eq:ras03} 
t_0=\frac{\xi-\cos^2\theta+\sqrt{\left(\cos^2\theta-\xi\right)^2+4\xi}}{2\xi},
\end{equation}
which satisfies $t_0>1$. Limiting values with respect to $\xi$ are
\begin{equation}\label{eq:ras04} 
\lim_{\xi\to0}t_0=\frac{1}{\cos^2\theta},\quad \lim_{\xi\to\infty}t_0=1.
\end{equation}

In the following analysis we assume that, next to large values of $r$, the parameter $\theta$ is well inside $(0,\frac12\pi)$, that is, $0<\delta\le \theta\le\frac12\pi-\delta$, where $\delta$ is  small and not depending on $r$. Also, we assume that $t_0\ge1+\delta$.

The integrand of the integral in \eqref{eq:ras01} has a pole at $t=t_p=1/y$, and we 
assume that  $t_0<t_p$. Later we will relax this assumption. Equality $t_0=t_p$ happens when $y=y_0$, where
\begin{equation}\label{eq:ras05} 
y_0=\frac{x+2p}{x+2p+2q}.
\end{equation}
Hence, we assume that $y<y_0$. This value $y_0$ is a transition value with respect to $y$: when $y\le y_0$ we have (roughly) $B_{p,q}(x,y)\le\bar{B}_{p,q}(x,y)$, and $B_{p,q}(x,y)$ is the smaller (primary) function for numerical computations; when 
$y_0\le y\le 1$, it is better to compute $\bar{B}_{p,q}(x,y)$ as the primary function, and use $B_{p,q}(x,y)=1-\bar{B}_{p,q}(x,y)$.

In terms of $x$ the transition value is
\begin{equation}\label{eq:ras06} 
x_0=2\frac{(p+q)y-p}{1-y}.
\end{equation}
$B_{p,q}(x,y)$ is the primary function for $x \ge x_0$.
For example, when we take $y=0.6$, $p=4.5$, $q=5.5$, we have $x_0=7.5$, and with $x=x_0$ we compute $B_{p,q}(x,y)\doteq0.4975$.

We use the substitution
\begin{equation}\label{eq:ras07} 
\phi(t)-\phi(t_0)=\tfrac12 w^2, 
\end{equation}
where we assume that for $t>1$ and $w\in\RR$: $\sign(t-t_0)=\sign(w)$. This transforms \eqref{eq:ras01} into
\begin{equation}\label{eq:ras08} 
B_{p,q}(x,y)=\frac{e^{-x/2+r\phi(t_0)}y^p (1-y)^q}{2\pi i}
\int_{-i\infty}^{i\infty}e^{\frac12rw^2} f(w)\,dw,
\end{equation}
where
\begin{equation}\label{eq:ras09} 
f(w)=\frac{1}{t(1-yt)}\,\frac{dt}{dw}, \quad \frac{dt}{dw}=\frac{w}{\phi^\prime(t)}.
\end{equation}

Substituting the expansion $\dsp{f(w)=\sum_{k=0}^\infty f_k w^k}$, we obtain
\begin{equation}\label{eq:ras10} 
B_{p,q}(x,y)\sim\frac{e^{-x/2+r\phi(t_0)}y^p (1-y)^q}{\sqrt{2\pi r}}
\sum_{k=0}^\infty (-1)^k f_{2k}\frac{2^k\left(\frac12\right)_k}{r^k}, \quad r\to\infty.
\end{equation}

This expansion can be used, uniformly with respect to $y\in[0,y_0-\delta]$, where $\delta$ is a small fixed positive number.
When the values of $y$, $p$, $q$ are given, the expansion can be used 
when $x > x_0$ (see \eqref{eq:ras06}).

The coefficients follow from the inversion of the transformation in \eqref{eq:ras07}. First we find the coefficients of the expansion
\begin{equation}\label{eq:ras11} 
t=t_0+t_1w+t_2w^2+\ldots,
\end{equation}
of which $t_1$ is given by 
\begin{equation}\label{eq:ras12} 
t_1=\left.\frac{dt}{dw}\right\vert_{w=0}=\frac{1}{\sqrt{\phi^{\prime\prime}(t_0)}}=\frac{t_0(t_0-1)}{\sqrt{\sin^2\theta\, t_0^2-(t_0-1)^2}}.
\end{equation}
The next ones are
\begin{equation}\label{eq:ras13} 
t_2=-\frac{\phi_3}{6\phi_2^2},\quad t_3=\frac{5\phi_3^2-3\phi_2\phi_4}{72\phi_2^{\frac72}},
\quad t_4=\frac{45\phi_4\phi_3\phi_2-40\phi_3^3-9\phi_5\phi_2^2}{1080\phi_2^5},
\end{equation}
where $\phi_k$ is the $k$th derivative of $\phi(t)$ at $t_0$.

The first coefficient $f_0$ of the expansion in \eqref{eq:ras10} is given by
\begin{equation}\label{eq:ras14} 
f_0=\frac{t_1}{t_0(1-yt_0)}=\frac{t_0-1}{(1-yt_0)\sqrt{\sin^2\theta\, t_0^2-(t_0-1)^2}}, 
\end{equation}
and the next one is
\begin{equation}\label{eq:ras15} 
\begin{array}{@{}r@{\;}c@{\;}l@{}}
f_2&=&\dsp{\frac{1}{24\phi_2^{\frac72}t_0^3(1-yt_0)^3}}\bigl(24\phi_2^2+12\phi_2(\phi_3-6y\phi_2)t_0\ +\\[8pt]
&&\quad (72\phi_2^2y^2+5\phi_3^2-36\phi_3y\phi_2-3\phi_4\phi_2)t_0^2\ +\\[8pt]
&&\quad 2y(12\phi_3y\phi_2+3\phi_4\phi_2-5\phi_3^2)t_0^3+y^2(5\phi_3^2-3\phi_4\phi_2)t_0^4\bigr).
\end{array}
\end{equation}

\subsubsection{The influence of the pole: error function approximation}\label{sec:erf}
When in \eqref{eq:ras01} the pole at $t=t_p=1/y$ is close to the saddle point $t_0$ the asymptotic approach has to be modified. For example, the shown coefficients $f_0$ and $f_2$  have poles when $t_0=t_p$. 

When we derived the expansion in \eqref{eq:ras10} we have assumed that $ t_0 < t_p$. This point $t_p$ has a corresponding point $w=\zeta$ that follows from the transformation in \eqref{eq:ras07}. When $t_p>t_0$, we have 
$\zeta>0$, when $t_p<t_0$ we should take $\zeta<0$. Hence,
\begin{equation}\label{eq:raserf01} 
\tfrac12 \zeta^2=\phi(t_p)-\phi(t_0),\quad \sign(\zeta)=\sign(t_p-t_0). 
\end{equation}

We write the function $f(w)$ defined in \eqref{eq:ras09} in the form
\begin{equation}\label{eq:raserf02} 
f(w)=\frac{A}{w-\zeta}+g(w),
\end{equation}
with the condition on $A$ that $g(w)$ is analytic at $w=\zeta$. We find that $A=1$ and because
\begin{equation}\label{eq:raserf03} 
e^{-r\phi(t_p)}= y^p(1-y)^q e^{-\frac12x}\quad \Longrightarrow\quad y^p(1-y)^q e^{-\frac12x}e^{r\phi(t_0)}=e^{-\frac12r\zeta^2},
\end{equation}
\eqref{eq:ras08} becomes
\begin{equation}\label{eq:raserf04} 
B_{p,q}(x,y)=\frac{e^{-\frac12r\zeta^2}}{2\pi i}
\int_{-i\infty}^{i\infty}e^{\frac12rw^2} \,\frac{dw}{w-\zeta} +
\frac{e^{-\frac12r\zeta^2}}{2\pi i}
\int_{-i\infty}^{i\infty}e^{\frac12rw^2} g(w)\,dw.
\end{equation}
Using \cite[7.2.3, 7.7.2]{Temme:2010:ERF}, we obtain
\begin{equation}\label{eq:raserf05} 
B_{p,q}(x,y)=\tfrac12\erfc\left(\zeta\sqrt{r/2}\right)+
\frac{e^{-\frac12r\zeta^2}}{2\pi i}
\int_{-i\infty}^{i\infty}e^{\frac12rw^2} g(w)\,dw,
\end{equation}
where $\erfc\,z$ is the complementary error function defined by
\begin{equation}\label{eq:raserf06} 
\erfc \,z=\frac{2}{\sqrt{\pi}}\int_z^\infty e^{-t^2}\,dt.
\end{equation}
For the complementary function introduced in \eqref{Bbar}  we find
\begin{equation}\label{eq:raserf07} 
\bar{B}_{p,q}(x,y)=\tfrac12\erfc\left(-\zeta\sqrt{r/2}\right)-
\frac{e^{-\frac12r\zeta^2}}{2\pi i}
\int_{-i\infty}^{i\infty}e^{\frac12rw^2} g(w)\,dw,
\end{equation}
where we have used $1-\frac12\erfc\,z=\frac12\erfc(-z)$.

An expansion of the integral now follows from substituting the expansion $\dsp{g(w)=\sum_{k=0}^\infty g_k w^k}$ to obtain
\begin{equation}\label{eq:raserf08} 
\begin{array}{@{}r@{\;}c@{\;}l@{}}
B_{p,q}(x,y)&=&\dsp{\tfrac12\erfc\left(\zeta\sqrt{r/2}\right) + 
\frac{e^{-\frac12r\zeta^2}}{\sqrt{2\pi r}}R_r(\zeta),} \\[8pt]
R_r(\zeta)&\sim&\dsp{\sum_{k=0}^\infty (-1)^k g_{2k}\frac{2^k\left(\frac12\right)_k}{r^k},\quad  r\to\infty.}
\end{array}
\end{equation}

We have derived this expansion under the condition $t_0<t_p$, that is, for $y<1/t_0$, but by taking $\zeta<0$ when $t_0>t_p$, we have (see \eqref{eq:raserf01}) $\sign(\zeta)=\sign(t_p-t_0)$, and we can use \eqref{eq:raserf08} for all $y\in(0,1)$.

For the asymptotic expansion in   \eqref{eq:raserf08} we assume that (see \eqref{eq:as01})
\begin{equation}\label{eq:raserf09} 
y\in[\delta_1,1-\delta_1], \quad \theta\in\left[\delta_2,\tfrac12\pi-\delta_2\right],
\quad \xi \ge 0,
\end{equation}
where $\delta_j$ are fixed small positive numbers.

When the argument of the complementary error function is large, we can use the asymptotic expansion
\begin{equation}\label{eq:raserf10} 
\tfrac12\erfc\left(\zeta\sqrt{r/2}\right) \sim 
\frac{e^{-\frac12r\zeta^2}}{\sqrt{2\pi r}}
\sum_{k=0}^\infty (-1)^k \frac{2^k\left(\frac12\right)_k}{r^k}\frac{1}{\zeta^{2k+1}},
\end{equation}
and we recover the expansion in \eqref{eq:ras10}. This means that the coefficients $g_{2k}$ follow from $f_{2k}$ by the relation
\begin{equation}\label{eq:raserf11} 
g_{2k}=f_{2k}-\frac{1}{\zeta^{2k+1}},\quad k=0,1,2,\ldots\,.
\end{equation}
When $\zeta\to0$ this relation remains valid, because the pole in $f_{2k}$ is properly cancelled analytically. The numerical evaluation of $g_{2k}$ for small values of $\vert\zeta\vert$ (that is, small values of $\vert t_p-t_0\vert$) needs some extra care.

\begin{remark}\label{rem:02}
{\rm
The variable $\zeta$ introduced in \eqref{eq:raserf01} is a function of $p$, $q$, $x$ and $y$. Clearly, $\zeta=0$ when saddle point and pole coincide:  $t_p=t_0$. This happens when $y$ and $x$ assume the transition values $y_0$ and $x_0$, see \eqref{eq:ras05}--\eqref{eq:ras06}. By using the simple relation
\begin{equation}\label{eq:raserf12} 
\frac{\partial \zeta}{\partial x}=\frac{y(t_p-t_0)}{2r\zeta},
\end{equation}
it is easy to derive the expansion
\begin{equation}\label{eq:raserf13} 
x=\sum_{k=0}^\infty x_k\zeta^k,
\end{equation}
where $x_0$ is the transition value $x_0=2(ry-p)/(1-y)$ (see \eqref{eq:ras06}). The next two coefficients are
\begin{equation}\label{eq:raserf14} 
x_1= \frac{2\sqrt{r(q-r(1-y)^2)}}{1-y},\quad
x_2= \frac{2(r(1-y)^3-q)}{3(r(1-y)^2-q)(1-y)}\,.
\end{equation}

The coefficients of the expansion  
\begin{equation}\label{eq:raserf15} 
y=\sum_{k=0}^\infty y_{k}\zeta^{k},
\end{equation}
can be obtained from
\begin{equation}\label{eq:raserf16} 
\frac{\partial \zeta}{\partial y}=\frac{xy(1-y)(t_p-t_0)+(2r+x)(y-y_0)}{2yr(1-y)\zeta}.
\end{equation}
Here, $y_0$ is the transition value $y_0=(x+2p)/(x+2r)$  (see \eqref{eq:ras05})  and the next  coefficients are
\begin{equation}\label{eq:raserf17} 
\begin{array}{@{}r@{\;}c@{\;}l@{}}
y_1&=& \dsp{-\frac{2\sqrt{q}\,\sqrt{4pr+4rx+x^2}}{\sqrt{r}\,(2r+x)^2},}\\[8pt]
y_2&=&\dsp{\frac{2(16pr^2(p-q)+8rx(8rp^2-p^2-3r^2)+12rx^2(p+r)+8rx^3+x^4)}{3r(2r+x)^3(4pr+4rx+x^2)}.}
\end{array}
\end{equation}
}
\end{remark}

\subsection{Computational examples}

As an example of the computational performance of
the expansions, Table~\ref{table2} shows relative errors in the computation of $B_{p,q}(x,y)$ using the asymptotic expansions
for the parameter values considered in Table 2 of \cite{Baharev:2008:CNF} and few additional values. 
When $q$ is a positive integer, we have in \eqref{eq:zas06} a finite sum. Then, as can be seen
in the Table, the result 
provided by the expansion is exact for $q=5$ by taking $n=4$.
In  \cite{Baharev:2008:CNF} also
a finite sum is given of a different form when their $b$ is a positive integer.

It is interesting to note that for few of the parameter values tested, the asymptotic expansion  
\eqref{eq:raserf08} provide an accuracy much better than single precision ($10^{-8}$)
just using three coefficients $g_{2k}$  \eqref{eq:raserf11}. 

 \begin{table}
$$
\begin{array}{lllllll}
\hline
 \mbox{AE} & p   & q & x & y & B^{AE}_{p,q}(x,y) & \mbox{Rel. error} \\ 
  \hline
\eqref{eq:zas06},\,n=5 & 2.3 & 3.5 & 54 & 0.8640 &   0.2760082728547706  & 5.6\,10^{-4} \\
\eqref{eq:zas06},\,n=5 & 2.3 & 3.5 & 140 & 0.9000 & 0.03608547984275312 & 2\,10^{-5}  \\
\eqref{eq:zas06},\,n=5  & 2.3 & 3.5 & 250 & 0.9000 & 0.0005034732632828640 & 8.9\, 10^{-7}  \\
\eqref{eq:zas06},\,n=4 & 5 & 5 & 54 & 0.8640 & 0.4563026193369792 & 0^{(*)} \\
\eqref{eq:zas06},\,n=4 & 5 & 5 & 140 & 0.9000 & 0.1041334930397555 & 0^{(*)}  \\
\eqref{eq:zas06},\,n=4  & 5 & 5 & 170 & 0.9560 & 0.6022421650011662 &0^{(*)}   \\
\eqref{eq:raserf08},\,k=2  & 10 & 10 & 54 & 0.8686 & 0.9187790583189610 & 5.7\,10^{-8} \\
\eqref{eq:raserf08},\,k=2 & 10 & 10 & 140 & 0.9000 &0.6008070986289955  & 1.4\,10^{-8}\\
\eqref{eq:raserf08},\,k=2 & 10 & 10 & 250 & 0.9000 & 0.09028986850391792 & 5.3\,10^{-7}\\
\eqref{eq:raserf08},\,k=2 & 20 & 20 & 54 & 0.8787 & 0.9998676573798253 & 9.0 \,10^{-12} \\
\eqref{eq:raserf08},\,k=2 & 20 & 20 & 140 & 0.9000 & 0.9925975041637949 & 5.2\,10^{-10}\\
\eqref{eq:raserf08},\,k=2 & 20 & 20 & 250 & 0.9220 & 0.9641190712607291  & 1.7\,10^{-9}\\
\eqref{eq:ras10},\,k=2 & 30 & 30 & 100 & 0.1 &  5.341313347397197\, 10^{-33}  & 3.7 \,10^{-6} \\
\eqref{eq:ras10},\,k=2 & 30 & 30 & 150 & 0.1 &   5.175358340461182\,10^{-42}   & 1.9\,10^{-6}\\
\eqref{eq:ras10},\,k=2 & 30 & 30 & 250 & 0.1 &    3.252685735589340 \,10^{-60}  & 7.8\,10^{-7}\\
\end{array}
$$
{\footnotesize {\bf Table~2}. Relative errors in the computation of $B_{p,q}(x,y)$ using the asymptotic expansions
\eqref{eq:zas06}, \eqref{eq:ras10} and \eqref{eq:raserf08} for different values of the parameters. The number of
terms used in each of the expansions is shown in the first column. 
$^{(*)}$ The asymptotic expansion is exact in this case. }
\label{table2}
\end{table}

\section{Inversion of the noncentral beta distribution}\label{sec:invert}

We consider the inversion problem of finding $x$ or $y$ when a value of $B_{p,q}(x,y)$ is given. 

As mentioned in \S\ref{sec:notations}, $x$ is the noncentrality parameter $\lambda$ of the noncentral $F$-distribution $F(w,\nu_1,\nu_2,\lambda)$. For $x=0$ we have $B_{p,q}(0,y)=I_y(p,q)$, the incomplete beta function. This follows from \eqref{defi}. Also, from \eqref{parx},  $B_{p,q}(x,y)$ is a decreasing function of $x$.
 Therefore when we consider the inversion for nonnegative values of $x$ with a given value $z$, we should write
\begin{equation}\label{eq:invert01} 
B_{p,q}(x,y)=z, \quad 0\le z\le I_y(p,q), \quad p, q, y, z \ {\rm given.}
\end{equation}

For the inversion with respect to $y$, we have the equation
\begin{equation}\label{eq:invert02} 
B_{p,q}(x,y)=z, \quad 0\le z\le 1, \quad  p, q, x, z \ {\rm given.}
\end{equation}
This corresponds with finding the quantiles of the distribution.

We try to find an approximation of $x$ or $y$ by using the representation in \eqref{eq:raserf08} and by inverting first the complementary error function, that is by finding $\zeta_0$ from the equation
\begin{equation}\label{eq:invert03} 
\tfrac12\erfc\left(\zeta_0\sqrt{r/2}\right)=z.
\end{equation}
A simple and efficient algorithm for computing the inverse of
the complementary error function is included, for example, in the package described
in \cite{Gil:2015:GCH}. Hence, we consider the numerical inversion  of (\ref{eq:invert03})   as a known problem 
and concentrate on finding from the computed value $\zeta_0$ the requested value of $x$ or $y$.

To obtain $x$ or $y$ when $\zeta_0$ is  small, we can use the expansions given in \eqref{eq:raserf13} or \eqref{eq:raserf15} with $\zeta=\zeta_0$, and for the general inversion problem we use \eqref{eq:raserf01}, written in the form
\begin{equation}\label{eq:invert04} 
\ln \left(\Frac{(2x)^{p/r}(S-2xy)^{q/r}}{(1-y)^{q/r}S}\right)+\Frac{2x-S}{4r}-\tfrac12\zeta_0^2=0,
\end{equation}
where  $S$ is defined by (see also \eqref{eq:ras03})
\begin{equation}\label{eq:invert05} 
t_0=\frac{S}{2xy}, \quad S=xy-2p+\sqrt{x^2y^2-4pxy+8rxy+4p^2}.
\end{equation}
After inverting equation \eqref{eq:invert03}, we obtain an approximation 
to the value of $x$ or $y$. 

A better approximation of $\zeta$ follows from obtaining $\zeta_1$ in the expansion $\zeta\sim\zeta_0+\zeta_1/r$. From \cite[Eq.~D.21]{Leeuwaarden:2009:AIE} we take\footnote{In the cited formula we should take $q=0$ for the present case.}
\begin{equation}\label{eq:invert06} 
\zeta_1=\frac{1}{\zeta_0}\ln\left(1+\zeta_0 g_0\right),
\end{equation}
where $g_0$  follows from  \eqref{eq:raserf11} and $f_0$ is given in \eqref{eq:ras14}.

\begin{figure}
\begin{center}
\epsfxsize=14cm \epsfbox{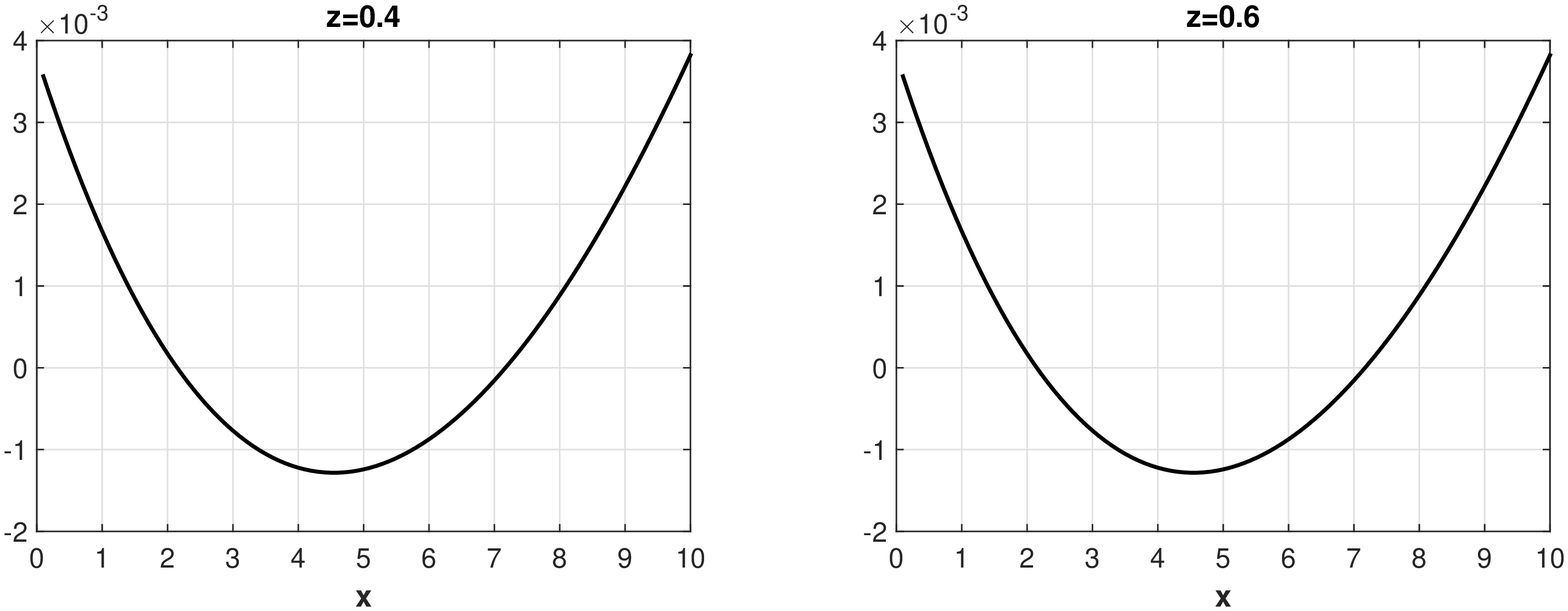}
\caption{ Plot of Eq.~\eqref{eq:invert04} written in the form $f(x)=0$, for
$p= 10$, $q=15$,  $y=0.45$ and $z=0.1,\,0.3$ ($z$ as in  \eqref{eq:invert01} and\eqref{eq:invert03}).
\label{fig:fig01}}
\end{center}
\end{figure}

\subsection{Examples for the inversion with respect to \protectbold{x}}\label{sec:invert01}
We take $p= 10$, $q=15$,  $y=0.45$. For these values the inversion problem makes sense if we take $z\le I_y(p,q) \doteq 0.7009$ (see \eqref{eq:invert01}). So, we can take $z=0.5$, which gives $\zeta_0=0$. The first approximation of $x$ that follows from the expansion in \eqref{eq:raserf13} is $x\doteq x_0=\frac{50}{11}\doteq4.54545$. With these values we have $B_{10,15}(x_0,0.45) \doteq 0.50952$.

Figure~\ref{fig:fig01} shows $f(x)$, the left-hand side of Eq.~\eqref{eq:invert04}, for two values of $z$ in Eq.~\eqref{eq:invert03} ($z^{(1)}=0.4$, which gives $\zeta_0^{(1)} \doteq 0.05067$, and $z^{(2)}=  0.6$, which gives  $\zeta_0^{(2)} \doteq -0.05067$).  The shown graphs in Figure~\ref{fig:fig01} attain their minimum value at the transition value $x_0=\frac{50}{11}$ (see \eqref{eq:ras06} and Remark~\ref{rem:02}). 

The function $f(x)$ shown in Figure~\ref{fig:fig01} has two zeros for both choices of $z^{(j)}$. For $z^{(1)}=0.4$, giving a positive $\zeta_0$, the zero that is larger than $x_0$ is an approximation of $x$,  because when $B_{p,q}(x,y)$ is smaller than $0.5$,  $x$ is larger than $x_0$ (see \eqref{eq:ras06}). For $z^{(2)}=0.6$ the corresponding $x$ is smaller than $x_0$.

The $x$-values resulting in the inversion of \eqref{eq:invert03} for these $z$-values
are $x^{(1)} \doteq 7.1704$, $x^{(2)} \doteq  2.1475$, respectively. In order to check the accuracy of the approximations,
we can compute the values of the noncentral beta function for $x^{(1)}$ and $x^{(2)}$: $B_{10,15}(x^{(1)},0.45) \doteq 0.40900$, 
$B_{10,15}(x^{(2)},0.45) \doteq 0.60915$.

The values $x^{(j)}$ for the parameters considered in the plots of Figure \ref{fig:fig01} are obtained by using the expansion of $x$  given in \eqref{eq:raserf13}. We have used terms up to $x_5\zeta_0^5$, and this term is smaller than $0.77e-5$. This shows the performance of the series.
As can be seen, in both cases it is possible to avoid the inversion of the nonlinear equation \eqref{eq:invert04} 
and obtain a good approximation to the $x$-values using  \eqref{eq:raserf13}. In general, these
approximations could also be used as starting values
of methods for solving nonlinear equations in case higher accuracy is required. 

The approximations of $x$ obtained in these examples have a relative error of about $0.03$. A better result can be obtained in examples with larger parameter $p$ and $q$; in the present case $r=p+q=25$, which is not a very large asymptotic parameter. Also,  in the examples only the complementary error function in representation \eqref{eq:raserf08} is used. A first influence of the $R_r$-term is contained in $\zeta_1$ given in \eqref{eq:invert06}. When we take this term into account, we find for the example with $z^{(1)}=0.4$, $\zeta^{(1)}\doteq0.05523$, giving $x^{(1)}\doteq7.4176$ and $B_{10,15}(x^{(1)},0.45) \doteq0.40014$, which is a better result.

\begin{figure}
\begin{center}
\epsfxsize=14cm \epsfbox{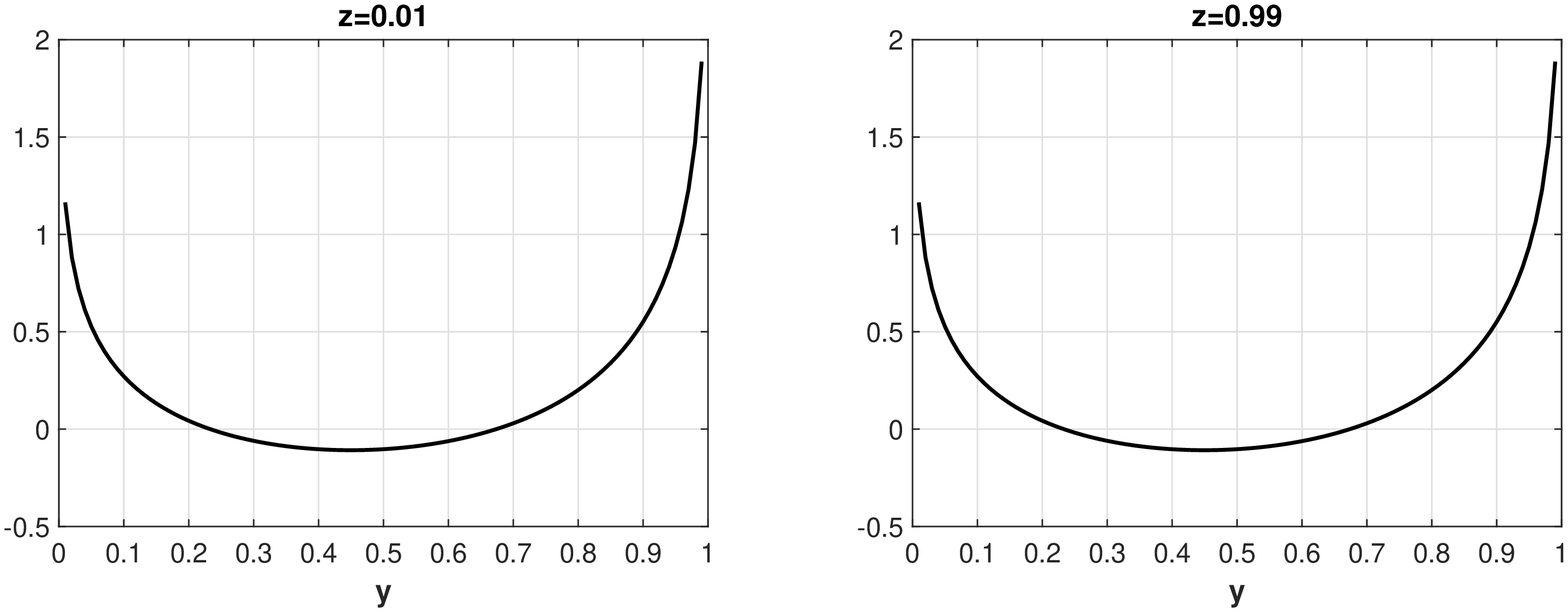}
\caption{ Plot of Eq.~\eqref{eq:invert04} written in the form $g(y)=0$, for
$p= 10$, $q=15$,  $x=4.5$ and $z=0.01,\,0.99$ ($z$ as in  \eqref{eq:invert01} and \eqref{eq:invert03}).
\label{fig:fig02}}
\end{center}
\end{figure}

\subsection{Examples for the inversion with respect to \protectbold{y}}\label{sec:invert02}
We take $p= 10$, $q=15$,  $x=4.5$. For $z=0.5$, the solution of \eqref{eq:invert03} is $\zeta_0=0$. The first approximation of $y$ that follows from the expansion in \eqref{eq:raserf15} is $y\doteq y_0= 49/109\doteq0.44954$. With these values we have $B_{10,15}(4.5,y_0) \doteq 0.5095$.

To compute $y$ from Eq.~\eqref{eq:invert04} when $x, p, q, z$ are given, we call the left-hand side of this equation $g(y)$. In Figure~\ref{fig:fig02} we show the graphs of $g(y)$ for values of $\zeta_0$ that follow 
from the inversion of Eq.~\eqref{eq:invert03} for $z^{(1)}=0.01$ and $z^{(2)}=0.99$. 

The values of $\zeta_0$ are 
$\zeta_0^{(1)} \doteq 0.4653$,  with $y^{(1)}\doteq0.2330$ (we take the zero left of $y_0=\frac{49}{109}\doteq 0.4495$ because, when $B_{p,q}(x,y)$ is smaller than $0.5$,  $y$ is smaller than $y_0$ (see \eqref{eq:ras05}). For $z^{(2)}=0.99$ the corresponding $y$ is larger than $y_0$. We have $\zeta_0^{(2)} \doteq -0.4652$,  with $y^{(2)}\doteq0.6739$. The computed values of the noncentral beta distribution are $B_{10,15}(4.5,y^{(1)}) \doteq0.011360$ and $B_{10,15}(4.5,y^{(2)}) \doteq0.98999$. 

We see a better result for the second $y$-value.  This can be explained by observing that, when the argument of  the complementary error function in \eqref{eq:raserf08} is negative and large in absolute value, this terms approaches 1 very fast, while the other term is exponentially small. In the present example the argument of  the complementary error function is $\zeta_0^{(2)} \sqrt{r/2}\doteq-1.6447$.

\section{Appendix}\label{sec:append}
We prove the relation in \eqref{eq:rel04} by using interchanging summation as in
\begin{equation}\label{eq:app01} 
\sum_{n=0}^\infty A_n \sum_{j=0}^n B_j=\sum_{j=0}^\infty B_j \sum_{n=j}^\infty A_n,
\end{equation}
where
\begin{equation}\label{eq:app02} 
A_n=  \frac{\Gamma(a+b+n)}{\Gamma(b)\Gamma(a+1+n)}x^n,\quad  B_j=\frac{\lambda^j}{j!}.
\end{equation}
We have
\begin{equation}\label{eq:app03} 
\begin{array}{@{}r@{\;}c@{\;}l@{}}
\dsp{\sum_{n=j}^\infty A_n}&=&\dsp{\sum_{n=0}^\infty A_{n+j}=x^j\sum_{n=0}^\infty \frac{\Gamma(a+b+n+j)}{\Gamma(b)\Gamma(a+1+n+j)}x^n}\\[8pt]
&=&\dsp{x^j\frac{\Gamma(a+b+j)}{\Gamma(b)\Gamma(a+1+j)}\sum_{n=0}^\infty \frac{(a+b+j)_n}{ (a+1+j)_n}x^n}\\[8pt]
&=&\dsp{\frac{x^j}{(a+j)B(a+j,b)}\FG{a+b+j}{1}{a+1+j}{x}}\\[8pt]
&=&x^{-a}(1-x)^{-b}I_x(a+j,b),
\end{array}
\end{equation}
where we have used \cite[8.17.8]{Paris:2010:IGA}. From this result, the relation in \eqref{eq:rel04} easily follows.

\section{Acknowledgements}
The authors thank the referees for their comments which have helped to improve the manuscript. 
We acknowledge financial support from 
{\emph{Ministerio de Ciencia e Innovaci\'on}}, project MTM2015-67142-P. 
NMT thanks CWI, Amsterdam, for scientific support.

\newpage

\bibliographystyle{plain}
\bibliography{biblio}

\end{document}